\documentclass{sigcomm-alternate}
\usepackage{amsfonts,latexsym,amsmath,amstext,amssymb,verbatim,epsfig,psfrag}
\usepackage{colordvi,pstcol}
\usepackage{graphicx,psboxit}
\usepackage{psfrag}

\def\R{\mbox{I\hspace{-.15em}R}}

\def\ol{\overline}

\newtheorem{remark}{Remark}
\newtheorem{theorem}{Theorem}

\begin{document}
\title{Optimizing the eigenvector computation algorithm with diffusion approach}

\numberofauthors{2}
\author{
   \alignauthor Dohy Hong\\
   \affaddr{Alcatel-Lucent Bell Labs}\\
   \affaddr{Route de Villejust}\\
   \affaddr{91620 Nozay, France}\\
   \email{\normalsize dohy.hong@alcatel-lucent.com}
   \alignauthor Philippe Jacquet\\
   \affaddr{Alcatel-Lucent Bell Labs}\\
   \affaddr{Route de Villejust}\\
   \affaddr{91620 Nozay, France}\\
   \email{\normalsize philippe.jacquet@alcatel-lucent.com}
}

\date{\today}
\maketitle

\begin{abstract}
In this paper, we apply the ideas of the matrix column based diffusion approach to
define a new eigenvector computation algorithm of a stationary probability
of a Markov chain.
\end{abstract}
\begin{psfrags}
\section{Introduction}\label{sec:intro}
In this paper, we assume that the readers are already familiar with the idea
of the fluid diffusion associated to the D-iteration \cite{d-algo}
to solve the equation:
$$
X = P.X + B
$$
and its application to PageRank equation \cite{dohy}.

For the general description of alternative or existing iteration methods, one
may refer to \cite{Golub1996, Saad}.

\section{Algorithm description}\label{sec:algo}
\subsection{Notation}
We recall that the D-iteration is defined by the couple $(P,B) \in \R^{N\times N}\times \R^N$ and exploits two state vectors:
$H_n$ (history) and $F_n$ (residual fluid):
\begin{eqnarray}
F_0 &=& B\\
F_n &=& (I_d - J_{i_n} + PJ_{i_n}) F_{n-1}.\label{eq:defF}
\end{eqnarray}
where $I_d$ is the identity matrix, $J_k$ a matrix with all entries equal to zero except for
the $k$-th diagonal term: $(J_k)_{kk} = 1$, $i_n$ the $n$-th node selected for the diffusion and
\begin{eqnarray}\label{eq:defH}
H_n &=& \sum_{k=1}^n J_{i_k} F_{k-1}\\
&=& \left(I_d - J_{i_n}(I_d - P)\right)H_{n-1} + J_{i_n} B.
\end{eqnarray}

The diffusion of a node $i$ containing $(F_n)_i = f$ means the following operation:
\begin{itemize}
\item $(H_n)_i += f$: this defined $H_{n+1}$;
\item $(F_n)_i = 0$ and $(F_n)_j += f\times p_{ji}$: this defines $F_{n+1}$.
\end{itemize}

Here, we'll use the equation satisfied by the state vectors (cf. \cite{d-algo}):
\begin{align}
H_n + F_n &= P. H_n + B.
\end{align}
We'll use here the notation $P.X$ meaning the usual matrix-vector product
$P\times X$.

We define $\sigma : \R^N \to \R$, by $\sigma(X)=\sum_{i=1}^N x_i$.
We define the $L_1$ norm, $|X|=\sum_{i=1}^N|x_i|$.

Below, we denote by $e$ the normalized unit column vector $1/N (1,..,1)^t$.

\subsection{Assumption}
We will assume in this paper that $P$ is a stochastic matrix (we took the
notation of the stochastic matrix per column), i.e.
$\sum_{i} p_{ij} = 1$ for each $i$. We define below the optimal algorithm
to find $X$ such that $X = P.X$.
For the sake of simplicity, we will assume that $P$ is irreducible and ergodic,
so that $P^n$ converges to a unique solution.
In fact, our approach finds a solution as soon as $P$ is stochastic and with no empty
columns or rows in $P$, or a bit more generally if $P^n . e$ is convergent.

\subsection{Algorithm description}
We first apply a full product $P.e$ then compensate by substracting $e$ so that
$\sigma(P.e-e) = 0$.
Then we set $F_0 = P.e - e$ and apply the diffusion iteration on $(P, F_0)$ (cf. \cite{d-algo}, \cite{revisit}).
This means that we solve $H$ such that $H = P.H + P.e-e$ and $X$ is obtained by $H + e$.
Note that $H+e$ is not necessarily normalized to 1.

\begin{remark}
The results we present here are in fact independent of the choice of $e$ and
we could use any other vector.
\end{remark}

\subsection{Some intuition on the optimality}
The benchmark test to existing methods will be addressed in a future paper.
From the intuition point of view, the D-iteration method was initially used
to solve $X = P.X + B$ when the spectral radius of $P$ is strictly less than
1, so that the fluid $B$ converges to zero exponentially.
For the problem $X=P.X$, the direct application of the diffusion would keep the
same amount of fluids in the system and the convergence can be obtained by Cesaro
averaging (or multiplicative normalization), but this would lead a very slow
convergence \cite{dohy}.
The algorithm above starts by creating an initial vector in the kernel of $\sigma$
in which we have positive fluids compensating exactly all negative fluids.
Then the convergence is obtained by the fact that positive fluids meeting negative
fluids just vanishes.

Now, intuitively, this is clearly optimal, because with the initial setting
we are {\em killing} all cycles and the diffusion behind {\em kills} all remaining forwarding
diffusions.

\section{Convergence}
\begin{theorem}
There exists a choice of the sequence of nodes such that
the diffusion applied on $(P, F_0)$ converges to $X-e$.
\end{theorem}
\proof
Let's first prove that the usual iteration method: $H_{n+1} = P.H_n + F_0$
is convergent. By induction, it is straightforward to obtain:
$H_n = P^n.e - e$.
Since $P^n.e\to X$, we have $H_{\infty} = H = X-e$.

For the diffusion process, from the equations \eqref{eq:defF} and \eqref{eq:defH}, we
can easily prove by induction that: $\sigma(F_n) = \sigma(H_n) = 0$. Then we prove that
$|F_n|$ is non-increasing function: when the diffusion is applied on a positive
fluid, the same amount is distributed; when they meet negative
fluid, a part of them vanishes so that $|F_n|$ is decreased, if not
$|F_n|$ is not modified.
For a more rigorous proof: set $j=i_{n+1}$ and$f=(F_n)_{i_{n+1}}$, then:
\begin{align*}
|F_{n+1}| =& \sum_{i\neq j} |(F_{n+1})_i| + |(F_{n+1})_j|\\
=& \sum_{i\neq j} |(F_n)_i + fp_{ij}| + |fp_{jj}|.
\end{align*}
Let's call $\Delta$ the set of nodes $i$ such that $(F_n)_i$ has a sign opposed to $f$.
Then,
\begin{align*}
|F_{n+1}| =& \sum_{i\neq j} (|(F_n)_i| + |f|p_{ij}) + |fp_{jj}|\\
    &+ \sum_{i\in \Delta} (|(F_n)_i + fp_{ij}| - |(F_n)_i| - |f|p_{ij})\\
=& |F_{n}| + \sum_{i\in \Delta} (|(F_n)_i + fp_{ij}| - |(F_n)_i| - |f|p_{ij}).
\end{align*}
Now, we use $|x+y| \le |x|+|y|$ to get $|F_{n+1}|\le |F_n|$.
Therefore, $|F_n|$ is convergent. 
The limit is necessarily equal to zero, because of the irreducibility of $P$:
there exists a path (diffusion sequence) such that
positive and negative fluids necessarily meet each others. 
For all $i,j$, there
exists $n$ such that $(P)^n_{ij} > 0$ implying there is a path $i_1=j,i_2,..,i_n=i$
with strictly positive weight $w(i,j) = p_{i_1,i_2}\times .. \times p_{i_{n-1},i_n} > 0$.
Applying the diffusion successively on the nodes $i_n,..i_2$ from $F_m=F$, we are sure to 
cancel at least $\min(|(F)_{i}|,w(i,j)\times |(F)_j|)$ (taking $(F)_i\times (F)_j < 0$).
Now taking $w = \min_{i,j} w(i,j)$, we are sure to cancel in less than $N$ diffusions
at least $w \times \max_i |(F)_i|$, which means an exponential convergence, and this
guarantee the convergence of $H_n$.


This theorem allows us to apply heuristic policy to optimize the sequence
of nodes for the diffusion.

\section{Application to PageRank equation}
The PageRank equation can be written as:
\begin{align}
X &= P.X
\end{align}
where $P = d\ol{P}_g + (1-d)/N J$, where $J$ is the matrix with all entries equal to 1
and $\ol{P}_g$ is the completed matrix from the initial $P_g$ (corresponding to the web graph)
by $e$ on columns associated to dangling nodes (nodes with no outgoing links).

If we apply the same method than above, we get:
$$
H = d\ol{P}_g H + d(\ol{P}_g.e - e).
$$
We can rewrite $d(\ol{P}_g.e - e)$ as $F_0 = d P_g .e - d\sigma(P_g.e) e$.

In this case, because of the factor $d$ and the presence of dangling nodes, the iteration
does not maintain $H_n, F_n$ in the kernel of $\sigma$.
However, the system is dominated by a system which is exponentially decreasing ($d$)
and it is easy to prove its convergence for D-iteration.

The limit of D-iteration on $(dP_g, F_0)$ satisfies:
$$
H = dP_g.H + F_0
$$
and 
\begin{align*}
H &= \sum_{i\ge 0} d^i P_g^i. F_0\\
  &= \sum_{i\ge 0} d^i P_g^i.[(d P_g.e - e) + (1-d +df_1) e]\\
  &= -e + \sum_{i\ge 0} d^i P_g^i.(1-d +df_1) e,
\end{align*}
where $f_1 = 1-\sigma(P_g.e)$.

The limit of D-iteration on $(d\ol{P}_g, F_0)$ satisfies:
\begin{align*}
H' &= d\ol{P}_g.H' + F_0 \\
   &= dP_g.H' + d \sigma(H'-P_g.H') e + F_0\\
   &= dP_g.H' + (d P_g.e - e) + (1-d +df_1+df_2) e\\
   &= dP_g.H' + (d P_g.e - e) + (1-d +df) e
\end{align*}
where $f_2 = \sigma(H'-P_g.H')$ and $f=f_1+f_2$,
and
\begin{align*}
H' &= \sum_{i\ge 0} d^i \ol{P}_g^i. F_0\\
  &= \sum_{i\ge 0} d^i P_g^i.[(d P_g.e - e) + (1-d +df) e]\\
  &= -e + \sum_{i\ge 0} d^i P_g^i.(1-d +df) e.
\end{align*}
Therefore, we have:
\begin{align*}
X &= H' + e\\
  &= \frac{1-d+df}{1-d+df_1}(H+e).
\end{align*}
Then $X$ can be computed from $H$ by normalizing $H+e$ to one.

\section{Extension to positive matrix}
Now, if $P$ is a positive irreducible matrix, Perron-Frobenius theorem says that
$(\rho^{-1} P)^n.e$ is convergent ($\rho$ is the spectral radius of $P$) to the unique
eigenvector which is strictly positive.
Therefore, to solve $P.X = \rho X$, we can apply the method above and solve
$H$ such that $H = P'.H + P.e - e$ with $P' = \rho^{-1} P$.
We would have $H_n = P'^n.e - e$ which converges to $X-e$. 

\begin{theorem}
There exists a choice of the sequence of nodes such that
the diffusion applied on $(P', P'.e-e)$ converges to $X-e$.
\end{theorem}
\proof
The proof is similar to the stochastic matrix case, except that we replace the $L_1$ norm
by the norm $|X|_V = \sum_i|x_i\times v_i|$ where $V$ is the left eigenvector of $P'$
(which has all entries strictly positive).
Then, $F_n$ is such that $\sigma_V(F_n) = \sum_i x_i\times v_i = 0$ and
$F_n$ is non-increasing function for the norm $|.|_V$.
Of course, the beauty of this result is that we don't need the explicit expression of
$V$.

\section{First convergence comparison}
\subsection{Error}
For the PageRank equation we will consider below, the distance to the limit
is computed as follows:
\begin{itemize}
\item for the power iteration (PI), the distance to the limit is
  bounded by $|X_{n+1}-X_n|\times d/(1-d)$; if $d=1$, we use an estimate
  of the distance to the limit by $|X_{n+1}-X_n|$;
\item for the initial D-iteration method (DI), the distance to the limit
 is given exactly by the $L_1$ norm of the residual fluid divided by $1-d$: $|F_n|/(1-d)$;
\item for DI$+$, the distance to the limit is bounded by $|F_n|/(1-d)$;
  for $d=1$, we use an estimate of the distance to the limit by  $|X_{n+1}-X_n|$.
\end{itemize} 
\subsection{Data set}
For the evaluation purpose, we used the 
web graph imported from the dataset \verb+uk-2007-05@1000000+
(available on \cite{webgraphit}) which has
41,247,159 links on 1,000,000 nodes.

Below we vary $N$ from $10^3$ to $10^6$ extracting from the dataset the
information on the first $N$ nodes.
Some graph properties are summarized in Table \ref{tab:1}:
\begin{itemize}
\item L: number of non-null entries (links) of $P$;
\item D: number of dangling nodes (0 out-degree nodes);
\item E: number of 0 in-degree nodes: the 0 in-degree nodes are defined recursively:
  a node $i$, having incoming links from nodes that are all 0 in-degree nodes, is
  also a 0 in-degree node; from the diffusion point of view, those nodes are those
  who converged exactly in finite steps;
\item O: number of loop nodes ($p_{ii} \neq 0$);
\item $\max_{in} = \max_i \#in_i$ (maximum in-degree, the in-degree of $i$ is the number of
  non-null entries of the $i$-th line vector of $P$);
\item $\max_{out} = \max_i \#out_i$ (maximum out-degree, the out-degree of $i$ is the number of
  non-null entries of the $i$-th column vector of $P$).
\end{itemize}

\begin{table}
\begin{center}
\begin{tabular}{|l|cccccc|}
\hline
N        & L/N  & D/N   & E/N   & O/N & $\max_{in}$ & $\max_{out}$\\
\hline
$10^3$   & 12.9 & 0.041 & 0.032 & 0.236 & 716   & 130\\
$10^4$   & 12.5 & 0.008 & 0.145 & 0.114 & 7982  & 751\\
$10^5$   & 31.4 & 0.027 & 0.016 & 0.175 & 34764 & 3782\\
$10^6$   & 41.2 & 0.046 & 0     & 0.204 & 403441& 4655\\
\hline
\end{tabular}\caption{Extracted graph: $N=10^3$ to $10^6$.}\label{tab:1}
\end{center}
\end{table}

To guarantee the continuity to $d=1$ and to meaningfully consider all nodes, the above graph has been completed
with one random outgoing link for all dangling nodes and with one random incoming link
for all nodes having no incoming links.

Note that for $d<1$, there exists a unique solution (stationary probability).
For $d=1$, DI is not defined, PI may not converge ($P$ may not be aperiodic) and
DI$+$ converges to the solution which is the limit of $d\to 1$.

\subsection{Comparison}
For the evaluation of the computation cost, we used
Linux (Ubuntu) machines: Intel(R) Core(TM)2 CPU, U7600, 1.20GHz, cache size 2048 KB (Linux1, $g++-4.4$).
The algorithms that we evaluated are:
\begin{itemize}
\item PI: Power iteration (equivalent to Jacobi iteration), using row vectors;
  for $d=1$, to force the convergence, we used relaxation idea with parameter $0.5$;
\item DI: D-iteration with node selection,
   if $(F_n)_i > r_{n'}\times \#out_i/L$, where $\#out_i$ is the out-degree of $i$
   and $r_{n'}$ is computed per cycle $n'$;
\item DI$+$: proposed solution. Initialization to $P.e-e$ followed by D-iteration with node selection,
   if $|(F_n)_i| > r_{n'}\times \#out_i/L$, where $\#out_i$ is the out-degree of $i$
   and $r_{n'}$ is computed per cycle $n'$.
\end{itemize}

\begin{table}
\begin{center}
\begin{tabular}{|l|cc|cc|}
\hline
         & nb iter & gain & time (s) & gain\\
\hline
\hline
 \multicolumn{5}{|l|}{$d=0.5$.} \\
\hline
PI     & 7    & x           & 0.01  & x\\
DI     & 4.8  & $\times$1.5 & 0     & $>\times$1\\
DI$+$  & 3.8  & $\times$1.8 & 0     & $>\times$1\\
\hline
 \multicolumn{5}{|l|}{$d=0.85$.} \\
\hline

PI     & 20 & x & 0.02  & x\\
DI     & 15.8  & $\times$1.3 & 0  & $>\times$2\\
DI$+$  & 11.8  & $\times$1.7 & 0  & $>\times$2\\
\hline
 \multicolumn{5}{|l|}{$d=0.99$.} \\
\hline
PI     & 469 & x & 0.14  & x\\
DI     & 202  & $\times$2.3 & 0.02  & $\times$7\\
DI$+$  & 193  & $\times$2.4 & 0.02  & $\times$7\\
\hline
 \multicolumn{5}{|l|}{$d=0.999$.} \\
\hline
PI & 5425 & x & 1.47  & x\\
DI     & 1379 & $\times$3.9 & 0.2  & $\times$7.4\\
DI$+$  & 1282 & $\times$4.2 & 0.17  & $\times$8.6\\
\hline
 \multicolumn{5}{|l|}{$d=1$.} \\
\hline
PI & 49& x &  0.03& x\\
DI  & undefined & x & undefined & x\\
DI$+$  & 242 & $\times$0.2 & 0.03  & $\times$1\\
\hline
\end{tabular}\caption{$N=10^3$: Comparison of the runtime for a target error of $1/N$. Gain: speed-up gain factor w.r.t. PI.}\label{tab:compa1}
\end{center}
\end{table}

\begin{table}
\begin{center}
\begin{tabular}{|l|cc|cc|}
\hline
         & nb iter & gain & time (s) & gain\\
\hline
\hline
 \multicolumn{5}{|l|}{$d=0.5$.} \\
\hline
PI & 9 & x & 0.04  & x\\
DI  & 5.3 & $\times$1.7 & 0.01  & $\times$4\\
DI$+$  & 5.1 & $\times$1.8 & 0  & $>\times$4\\
\hline
 \multicolumn{5}{|l|}{$d=0.85$.} \\
\hline
PI & 42 & x & 0.15  & x\\
DI  & 15 & $\times$2.8 & 0.02  & $\times$8\\
DI$+$  & 13.6 & $\times$3.1 & 0.01  & $\times$15\\
\hline
 \multicolumn{5}{|l|}{$d=0.99$.} \\
\hline
PI & 541 & x & 1.73  & x\\
DI  & 120 & $\times$4.5 & 0.16  & $\times$11\\
DI$+$  & 117 & $\times$4.6 & 0.15  & $\times$12\\
\hline
 \multicolumn{5}{|l|}{$d=0.999$.} \\
\hline
PI & 7739 & x & 24.4  & x\\
DI  & 919 & $\times$8.4 & 1.34  & $\times$18\\
DI$+$  & 880 & $\times$8.8 & 1.25  & $\times$20\\
\hline
 \multicolumn{5}{|l|}{$d=1$.} \\
\hline
PI & 393 & x & 1.26  & x\\
DI  & undefined & x & undefined & x\\
DI$+$  & 53 & $\times$7.4 & 0.07  & $\times$18\\
\hline
\end{tabular}\caption{$N=10^4$: Comparison of the runtime for a target error of $1/N$. Gain: speed-up gain factor w.r.t. PI.}\label{tab:compa2}
\end{center}
\end{table}

\begin{table}
\begin{center}
\begin{tabular}{|l|cc|cc|}
\hline
         & nb iter & gain & time (s) & gain\\
\hline
\hline
 \multicolumn{5}{|l|}{$d=0.5$.} \\
\hline
PI & 12 & x & 1.1  & x\\
DI     & 5.7 & $\times$2.1 & 0.14  & $\times$8\\
DI$+$  & 5.9 & $\times$2.0 & 0.14  & $\times$8\\
\hline
 \multicolumn{5}{|l|}{$d=0.85$.} \\
\hline
PI & 53 & x & 4.6  & x\\
DI     & 17.4 & $\times$3.0 & 0.44  & $\times$10\\
DI$+$  & 16.6 & $\times$3.2 & 0.41  & $\times$11\\
\hline
 \multicolumn{5}{|l|}{$d=0.99$.} \\
\hline
PI & 834 & x & 72  & x\\
DI & 163 &  $\times$5.1 & 4.6 &$\times$16\\
DI$+$  & 156 & $\times$5.3 & 4.3  & $\times$17\\
\hline
 \multicolumn{5}{|l|}{$d=0.999$.} \\
\hline
PI & 8253 & x & 710  & x\\
DI     & 936 & $\times$8.0 & 32  & $\times$22\\
DI$+$  & 897 & $\times$9.2 & 29  & $\times$24\\
\hline
 \multicolumn{5}{|l|}{$d=1$.} \\
\hline
PI & 15546 & x & 1357 & x\\
DI  & undefined & x & undefined & x\\
DI$+$  & 931 & $\times$17 & 30  & $\times$45\\
\hline
\end{tabular}\caption{$N=10^5$: Comparison of the runtime for a target error of $1/N$. Gain: speed-up gain factor w.r.t. PI.}\label{tab:compa3}
\end{center}
\end{table}

\begin{table}
\begin{center}
\begin{tabular}{|l|cc|cc|}
\hline
         & nb iter & gain & time (s) & gain\\
\hline
\hline
 \multicolumn{5}{|l|}{$d=0.5$.} \\
\hline
PI    & 15   & x  & 17.1 & x\\
DI    & 6.7  & $\times$2.2 & 1.9 & $\times$9\\
DI$+$ & 6.6  & $\times$2.2 & 1.9  & $\times$9\\
\hline
 \multicolumn{5}{|l|}{$d=0.85$.} \\
\hline
 PI  & 65   & x & 74 & x  \\
 DI    & 18.1  &$\times$3.6  & 5.8 & $\times$13\\
 DI$+$ & 18.0  &$\times$3.6  & 5.5 & $\times$13\\
\hline
 \multicolumn{5}{|l|}{$d=0.99$. } \\
\hline
 PI    & 1027  & x          & 1167 & x \\
 DI    & 171  & $\times$6.0 & 66 & $\times$18  \\
 DI$+$ & 168  & $\times$6.1 & 60 & $\times$19  \\
\hline
 \multicolumn{5}{|l|}{$d=0.999$.} \\
\hline
 PI  & 10783 & x  & 12235   & x \\
 DI    &  991  & $\times$11 & 455 & $\times$27  \\
 DI$+$ &  984  & $\times$11 & 433 & $\times$28 \\
\hline
 \multicolumn{5}{|l|}{$d=1$.} \\
\hline
 PI  & 105700 & x  &  121000  & x \\
 DI  & undefined & x & undefined & x\\
 DI$+$ & 3087   & $\times$34 & 1605 & $\times$75  \\
\hline
\end{tabular}\caption{$N=10^6$: Comparison of the runtime for a target error of $1/N$. Gain: speed-up gain factor w.r.t. PI.}\label{tab:compa4}
\end{center}
\end{table}


\section{Conclusion}\label{sec:conclusion}
We proposed a new algorithm to accelerate the computation of the dominant eigenvector of non-negative
matrix inspired from the D-iteration diffusion vision.
To show its potential, first evaluation results are included.

\end{psfrags}
\bibliographystyle{abbrv}
\bibliography{sigproc}

\end{document}